\documentclass[11pt]{article}
\usepackage[T1]{fontenc}
\usepackage[latin9]{inputenc}
\usepackage[a4paper]{geometry}
\usepackage{color}
\usepackage{amsmath}
\usepackage{setspace}
\usepackage{amssymb}

\usepackage[active]{srcltx}

\newtheorem{definition}{Definition}[section]
\newtheorem{remark}[definition]{Remark}
\newtheorem{theorem}[definition]{Theorem}

\newtheorem{proposition}{Proposition}[section]

\newcommand{\modifr}[1]{\textcolor{red}{}}

\def\Fc{\stackrel{\,\,\circ}{F}}
\def\Fb{\stackrel{\,\,\bullet}{F}}

\def\Dc{\stackrel{\,\,\circ}{\Delta}}
\def\Db{\stackrel{\,\,\bullet}{\Delta}}

\newcommand{\R}{\mathbb{R}}

\newcommand{\Z}{\mathbb{Z}}
\newcommand{\E}{\mathrm{E}}
\newcommand{\conj}[1]{\overline{#1}}

\begin{document}

\title{\textbf{\Large On a new resolution of the identity for coherent states of the Euclidean group $E(2)$ }}

\author{Julio Guerrero$^{a,b}$}

\maketitle
\begin{singlespace}

{\noindent \begin{center}
\emph{$^{a}$ 
Department of Mathematics, Faculty of  Experimental Sciences, University of Ja\'en, Campus Las Lagunillas s/n, 23071 Ja\'en, Spain}
\end{center}}

\vspace{-0.75cm}

{\noindent \begin{center}
\emph{$^{b}$ 
Institute Carlos I of Theoretical and Computational Physics, University of  Granada, Fuentenueva s/n, 18071 Granada, Spain}
\end{center}}

\end{singlespace}

\begin{abstract}
We provide a non-divergent resolution of the identity for coherent states of the Euclidean group $E(2)$, replacing the integral on the whole group with the Haar measure by an invariant integration
with non local kernels constructed by Steinberg and Wolf in 1981, valid for solutions of Helmholtz equation, that appears in natural way as the eigenvalue equation of the Casimir operator of $E(2)$.
With this construction, all properties of Perelomov coherent states are preserved, in particular covariance, and all square-integrable functions are addmissible, improving previous results like those of de Bi\`evre and Torr\'esani.

\end{abstract}
PACS: \\
Keywords:  Coherent states, Euclidean group, Helmholtz equation


\section{Introduction}
\label{introduction}

 Let $G$ be a locally compact group, with $d\mu^L(g)$ its (left) Haar measure.  Let $U$ be a unitary and irreducible representation of $G$ on a Hilbert space $\cal H$. In the following,  $U$ will be fixed, and we shall introduce functions $\eta\in\cal H$ to construct certain functions and operators on $\cal H$ and $L^2(G)$ that will depend on $U$ and on $\eta$. To simplify the notation, we shall label them with $\eta$, omiting the dependence on $U$. Also, as otherwise stated, the scalar produt $\langle\cdot,\cdot\rangle$ with be the scalar product in ${\cal H}$.
 
 The notion of coherent states of Gilmore-Perelomov type is intimately linked to the notion of square integrability of unitary representations and to the existence of admissible functions.

 \begin{definition}
 A function $\eta\in\cal H$ is denoted admissible and the representation $U$ square integrable if
 \begin{equation}
  c_\eta=\int_{G}d\mu^L(g) |\langle U(g)\eta,\eta\rangle|^2<\infty \,.
 \end{equation}
\end{definition}

\begin{definition}
Let $ \eta\in {\cal H}$ admissible, the coherent states associated to $G$ related to $ \eta$ are functions in ${\cal H}$ labeled by the elements of $G$  given by:
\begin{equation}
\eta_{g}=  U(g) \eta\,\qquad g\in G\,.
\end{equation}
\label{defCS}
\end{definition}

In the case when there exists admissible functions, the CS map is well defined.
\begin{proposition} \label{CSmap}
Suppose $\eta\in {\cal H}$ is admissible. Then the CS map given by
 \begin{equation}
  f\in{\cal H}\rightarrow F(g)\equiv T_\eta f(g)= \langle \eta_{g},f\rangle  
 \end{equation}
 is an application $T_\eta:{\cal H}\rightarrow L_\eta^2(G)\equiv  {\rm Range}(T_\eta)  \subset  L^2(G)$.

\end{proposition}

See \cite{Grossmann85} for a proof. The subspace  $L_\eta^2(G)$ is a reproducing kernel Hilbert space (RLHS), see Proposition \ref{OverlappingKernel}.

The family of coherent states $\eta_g$ is covariant under the action of the group $G$, $U(g')\eta_g=\eta_{g'g}$, as derived immediately from Definition \ref{defCS},  and this allows to prove  the following proposition.

\begin{proposition}\label{intenrwinT}
The coherent state map interwines the representation $U$ and the left regular representation $\Pi^L$ on
$L^2(G)$:
\begin{equation}
   \Pi^L(g) T_\eta = T_\eta U(g)\,,\qquad \forall g\in G
\end{equation}
\end{proposition}

See \cite{Grossmann85} for a proof. This proposition implies also that $L_\eta^2(G)$ is isomorphic to a subspace of $L^2(G)$
carrying an irreducible subrepresentation of $\Pi^L$, that we shall denote by $\Pi^L_\eta$, which is unitarily equivalent to $U$. This proposition will be crucial for the main result of the paper.

\begin{proposition}\label{adjointCSmap}
The adjoint of $T_\eta$,  $T_\eta^\dagger\,: L^2(G)\rightarrow {\cal H}$, is given by:
\begin{equation}
 T_\eta^\dagger \psi =\,\int_{G}d\mu^L(g) \psi(g) \eta_g
\end{equation}

\end{proposition}

The adjoint  $T_\eta^\dagger$ verifies a similar property, due to the invariance of the left Haar measure on $G$:

\begin{proposition} \label{intenrwinTd}
 The adjoint $T_\eta^\dagger\,: L^2(G)\rightarrow {\cal H}$ also interwinnes $U$ and $\Pi^L$:
 \begin{equation}
 T_\eta^\dagger \Pi^L(g)=U(g) T_\eta^\dagger 
 \end{equation}

\end{proposition}

Also, it is possible to provide a resolution of the identity.
\begin{proposition} \label{Abounded}
 Suppose $\eta\in {\cal H}$ is admissible. Then the operator $A_\eta=T_\eta^\dagger T_\eta$ in $\cal H$ given by 
 \begin{equation}
  A_\eta f=\int_{G}d\mu^L(g) \langle \eta_g,f\rangle \eta_g
 \end{equation}
 is bounded, invertible and with a bounded inverse.

\end{proposition}

See \cite{AAG14} for a proof. If $U$ is irreducible, $A$ is proportional to the Identity operator in $\cal H$ if $\eta$ is admissible:

\begin{theorem}\label{Aidentity}
 Suppose $\eta\in {\cal H}$ is admissible and $U$ irreducible. Then the operator $A_\eta=T_\eta^\dagger T_\eta=\lambda_\eta I_{\cal H}$, with $0<\lambda_\eta<\infty$.
\end{theorem}

\noindent \textit{Proof:} (see \cite{Grossmann85}) From Prop. \ref{intenrwinT} and \ref{intenrwinTd}, we have:
\begin{equation}
 A_\eta U(g) = T_\eta^\dagger T_\eta U(g) = T_\eta^\dagger \Pi^L(g) T_\eta = U(g) T_\eta^\dagger T_\eta = U(g) A_\eta
\end{equation}
and therefore $A$ interwinnes the irreducible representation $U$. By Schur's Lemma $A_\eta=\lambda_\eta I_{\cal H}$. Since $\eta$ is admissible, by Prop.  \ref{Abounded} $A_\eta$ is bounded with a bounded inverse and therefore $0<\lambda_\eta<\infty$.$\square$

\begin{remark}\label{Remark-ResolutionIdentity}
 Note that for a tight frame we have a resolution of the identity:
 \begin{equation}
  \langle h, \frac{1}{\lambda_\eta}A_\eta f\rangle=\frac{1}{\lambda_\eta}\int_{G}d\mu^L(g) \langle \eta_g,f\rangle\langle h,\eta_g \rangle = \frac{1}{\lambda_\eta}\int_{G}d\mu^L(g)    \langle \eta_g,f\rangle\conj{\langle \eta_g, h \rangle}=\langle f,h\rangle
 \end{equation}

\end{remark}

It is possible to construct another important operator from $T_\eta$ and $T_\eta^\dagger$.

\begin{proposition}\label{OverlappingKernel}
 Suppose $\eta\in {\cal H}$ is admissible. Define the  operator $B_\eta=T_\eta T_\eta^\dagger\,:\,L^2(G)\rightarrow L_\eta^2(G)$. Then $B_\eta$ is a convolution operator which is bounded, hermitian,  positive definite on  $L_\eta^2(G)$ and with integral kernel $K_\eta(g,g')=K_\eta(g'^{-1}g)$ satisfying the reproducing property:
 \begin{equation}
  \int_G d\mu^L(g') K_\eta(g'^{-1}g)\Psi(g') = \Psi(g) \,,\forall \Psi\in   L_\eta^2(G)
 \end{equation}

\end{proposition}

\noindent\textit{Proof:} See \cite{AAG14}.

\begin{remark}\label{Remark1-OverlappingKernel}
Note that, from the definition of $T_\eta$ and $T_\eta^\dagger$, we have:
\begin{equation}
 K_\eta(g'^{-1}g) = \langle \eta_g,\eta_{g'}\rangle
\end{equation}
and for this reason it is usually known as the overlapping kernel.

\end{remark}

\begin{remark}\label{Remark2-OverlappingKernel}
Note that $B_\eta$, if suitably normalized, is an orthogonal projector from $L^2(G)$ onto $L_\eta^2(G)$. And the kernel  $K_\eta$ plays the role of a Dirac delta on $L_\eta^2(G)$.

\end{remark}

\begin{remark}\label{Remark3-OverlappingKernel}
Note that $B_\eta$ endowes  $L_\eta^2(G)$ with a Reproducing Kernel Hilbert (RKHS) space structure. 

\end{remark}

\section{The Euclidean group $E(2)$ and its unitary and irreducible representations}

Let's consider the usual unitary and irreducible representations of the Euclidean group $E(2)$ on the Hilbert space ${\cal H}=L^2([0,2\pi))$ of square integrable functions on the circle \cite{Isham91}, considered as $2\pi$-periodic functions on $\mathbb{R}$, that can be characterized by their values on $[0,2\pi)$:
\begin{equation}
 U_{(\varphi,v,w)}^\lambda f(\theta)=e^{i \lambda(v\sin\theta + w\cos\theta)}f(\theta-\varphi)\,,\qquad \lambda>0\,.
 \end{equation}
 
 Note that in the literature there is another (unitarily equivalent) representation given by:
 \begin{eqnarray}
  \tilde U_{(\varphi,v,w)}^\lambda f(\theta) &=& e^{i \lambda(v\sin(\theta-\varphi) + w\cos(\theta-\varphi))}f(\theta-\varphi)\nonumber \\ &=& U^\lambda_{(\varphi,v\cos\varphi+w\sin\varphi,w\cos\varphi-v\sin\varphi)}f(\theta)   \\
  &=& U^\lambda_{(\varphi,v',w')} f(\theta) \nonumber
 \end{eqnarray}
 where 
 \begin{equation}
  \left(\begin{array}{l}
         v'\\
         w'
        \end{array}\right)= R(\varphi) \left(\begin{array}{l}
         v\\
         w
        \end{array}\right)
\end{equation}
with $R(\varphi)$ a rotation matrix of angle $\varphi$.  We shall use the representation $U^\lambda$ for simplicity, but all the results obtained are also valid for the representation  $\tilde U^\lambda$.

\section{Coherent states for $E(2)$}

We shall introduce coherent states (CS) for $E(2)$ in the usual way, as the action of the unitary operators $U^\lambda$ on a fixed function in ${\cal H}$.

\begin{definition}
Let $ \eta\in {\cal H}$, the coherent states associated to $E(2)$ related to $ \eta$ are given by:
\begin{equation}
\eta_{(\varphi,v,w)}(\theta)=  U^\lambda_{(\varphi,v,w)} \eta(\theta)= e^{ i \lambda(v\sin\theta + w\cos\theta)}\eta(\theta-\varphi) \label{CSE2}
\end{equation}

\end{definition}

In order to have a useful family of coherent states, it is necessary that the fixed function $\eta$ be admissible. However, it is well-known \cite{Isham91}  that there are no admissible functions for $E(2)$ (neither for the Euclidean groups $E(n)\,,n\geq 2$).
\begin{proposition} \label{NonAdmissible}
 For any $ \eta\in {\cal H}$,
 \begin{equation}
  \int_{E(2)}d\mu^L(\varphi,v,w) |\langle \eta_{(\varphi,v,w)},\eta\rangle|^2=+\infty
 \end{equation}
 where $d\mu^L(\varphi,v,w)=d\varphi dvdw$ is the left Haar measure for $E(2)$.

\end{proposition}
\textit{Proof:} The proof of this proposition can be found in Appendix \ref{app:0}.

The non-existence of admissible functions is related on the divergence on the integral with respect to both  $v$ and $w$, which diverge for any function $\eta\in{\cal H}$. The unique hope for the construction of coherent states for $E(2)$ with admissible functions is to restrict the integration to some subset of coherent states which depend on only one of the variables $v$ or $w$ (or a mix of them).

In this line goes the proposal of de Bi\`evre \cite{DeBievre89}  (see also  \cite{Torresani95}),
where the coherent states are restricted to depend only on the \textit{tangential momentum} (the variables $v$ and $w$ can be considered as momentum variables of the angle variable $\varphi$). For instance, at the north pole $(0,1)$ of the circle (identified as the circunferencie or radiuos 1 in the plane), the tangential momentum is $(v,0)$ (while $(0,w)$ is normal to the circle). Of course, the tangential and normal momenta change from point to point of the circle.

In \cite{Gazeau18} the proposal of de Bi\`evre is formalized in the setting of coherent states modulo a subgroup \cite{AAG14}.

The main problem with de Bi\`evre's proposal is that the addmissible function $\eta$ must be restricted to have support on $(-\frac{\pi}{2},\frac{\pi}{2})$ (the north hemisphere), and the admissibility condition reads (for a proof see \cite{DeBievre89,Torresani95,Gazeau18}):
\begin{proposition}
 A function $\eta\in{\cal H}$ with $supp(\eta)\subset [-\frac{\pi}{2}, \frac{\pi}{2}]$ such that 
 \begin{equation}
  0<c_\eta=\int_{-\pi/2}^{\pi/2}d\theta  \frac{|\eta(\theta)|^2}{\cos\theta} <\infty
 \end{equation}
 is admissible when the family of coherent states is restricted to $\eta_{(\varphi,v,0)}$. Furthermore, the family of coheret states provide a resolution of the identity:
 \begin{equation}
 \frac{1}{c_\eta} \int_{-\pi}^\pi d\varphi \int_{-\infty}^\infty dv \langle f,\eta_{(\varphi,v,0)}\rangle
  \langle \eta_{(\varphi,v,0)},g\rangle = \langle f,g\rangle\,,\qquad \forall f,g\in{\cal H}
 \end{equation}

 \end{proposition}
 
 We shall introduce a new family of coherent states, which are not restricted to depend on any subset of $E(2)$, but  with a completely different resolution of the identity.
The main motivation for the introduction of our family of coherent states is the following rather natural, but susprisingly unknown result.
 
 \begin{theorem}\label{Helmholtz}
 Let $\eta\in{\cal H}$, then the coherent states $\eta_{(\varphi,v,w)}$ verify Helmholtz equation:
 \begin{equation}
  \left(\frac{\partial^2\,}{\partial v^2}+\frac{\partial^2\,}{\partial w^2}+\lambda^2\right)\eta_{(\varphi,v,w)}=0 \label{EqHelmholtz}
 \end{equation}
 \end{theorem}
 
 \textit{Proof:}
Consider the following vector fields on $\R^2$:
\[
  T_v= \frac{\partial}{\partial v},\qquad
  T_w= \frac{\partial}{\partial w},\qquad
  J=\left(v\frac{\partial}{\partial w}-w\frac{\partial}{\partial v}\right).
\]

They close the Lie algebra of $E(2)$:
\[
  [T_v,T_w]=0,\qquad
  [J,T_v]=-T_w,\qquad
  [J,T_w]=T_v.
\]

These operators are the infinitesimal generators of the representation $U^\lambda$ restricted to $E(2)/SO(2)$.

Given $f(\theta)\in L^2([0,2\pi])$, we define
\[
  f(\varphi,v,w)=\langle\eta_{(\varphi,v,w)},f\rangle
  =\int_{-\pi}^{\pi}d\theta\;e^{-i\lambda(v\sin\theta+w\cos\theta)}\,
   \conj{\eta(\theta-\varphi)}\,f(\theta).
\]

Then

\begin{eqnarray}
  T_v\,f(\varphi,v,w)
  &=\frac{\partial}{\partial v}\langle\eta_{(\varphi,v,w)},f\rangle
   =\frac{\partial}{\partial v}\int_{-\pi}^{\pi}d\theta\;
    e^{-i\lambda(v\sin\theta+w\cos\theta)}\,\conj{\eta(\theta-\varphi)}f(\theta)\nonumber\\
  &=-i\lambda\int_{-\pi}^{\pi}d\theta\;\sin\theta\;
    e^{-i\lambda(v\sin\theta+w\cos\theta)}\,\conj{\eta(\theta-\varphi)}f(\theta)\\[4pt]
  T_w\,f(\varphi,v,w)
  &=\frac{\partial}{\partial w}\langle\eta_{(\varphi,v,w)},f\rangle
   =\frac{\partial}{\partial w}\int_{-\pi}^{\pi}d\theta\;
    e^{-i\lambda(v\sin\theta+w\cos\theta)}\,\conj{\eta(\theta-\varphi)}f(\theta)\nonumber\\
  &=-i\lambda\int_{-\pi}^{\pi}d\theta\;\cos\theta\;
    e^{-i\lambda(v\sin\theta+w\cos\theta)}\,\conj{\eta(\theta-\varphi)}f(\theta)
\end{eqnarray}

Therefore
\begin{eqnarray}
  \left(T_v^2+T_w^2\right)f(\varphi,v,w)
  &=-\lambda^2\int_{-\pi}^{\pi}d\theta\,(\sin^2\theta+\cos^2\theta)\;
    e^{-i(v\sin\theta+w\cos\theta)}\,\conj{\eta(\theta-\varphi)}f(\theta)\nonumber\\
  &=-\lambda^2f(\varphi,v,w).
\end{eqnarray}
Thus
\[
  \left(\partial_v^2+\partial_w^2+\lambda^2\right)f(\varphi,v,w)=0\,.
\]
We conclude that $f(\varphi,v,w)$ is a solution of the Helmholtz equation with
respect to $(v,w)$ with \textit{wave number} $\lambda$. Note that $\eta_{(\varphi,v,w)}(\theta)$ also satisfies this property with respect to $(v,w)$.
 

 This property is inherited by the coherent state map $T_\eta$ introduced in Prop. \ref{CSmap}, in such a way that the functions
  $F(\varphi,v,w)=\langle \eta_{(\varphi,v,w)},f\rangle$ on $E(2)$ also verify Helmholtz equation, although they don't belong to $L^2(G)$ by Proposition \ref{NonAdmissible}.

This theorem has very far reaching consecuences, as we shall see along this paper. Theorem \ref{Helmholtz} can be seen as a natural consequence\footnote{Althoug the representation $U^\lambda$ is not square integrable, the Prop. \ref{intenrwinT} remains valid if a suitable scalar product is introduced in $Range(T_\eta)$, see later.}  of Prop. \ref{intenrwinT} and the irreducibilty of the representation $U^\lambda$. In fact $C_2= \frac{\partial^2\,}{\partial v^2}+\frac{\partial^2\,}{\partial w^2}$ is the Casimir of the Euclidean group $E(2)$ for the left (and right) regular representation, and Helmholtz equation is nothing other than the eigenvalue equation   $C_2=\lambda^2 I_{L^2(E(2))}$.
 
 This eigenvalue equation simply states that the irreducible representation $U^\lambda$ is related (in the sense of the Coadjoint Orbit method \cite{Kirillov76}) to the particular coadjoint orbit
 given by a cylinder of radius $\lambda$ in $\mathbb{R}^3$.

 Functions in the Hilbert space ${\cal H}$ carrying the representation $U^\lambda$ implicitly satisfy this eigenvalue equation, since they are already defined on the circle. However, the functions  $F(\varphi,v,w)$  on $E(2)$ must  satisfy the Casimir eigenvalue equation with eigenvalue $\lambda^2$ in order to belong\footnote{Remember that, to avoid cluttering the notation, we did not specify the dependence on $U^\lambda$ in all objects constructed in Sec. \ref{introduction}, and therefore $T_\eta$ also depends on $\lambda$.} to $Range(T_\eta)$.
 
 \begin{definition}
  Let's denote by ${\cal F}_\lambda$ the vector space of oscillatory solutions of Helmholtz equation (\ref{EqHelmholtz}). Oscillatory solutions of Helmholtz equation are those whose 2D Fourier transform have support on the disk of radius $\lambda$. We will consider in ${\cal F}_\lambda$ functions in $E(2)$ such that are oscillatory solutions of Helmholtz equation with respect to the variables $(v,w)$.
 \end{definition}

 For a detailed account on oscillatory solutions of Helmholtz equation, see Appendix B of \cite{SphereMomentum20}, where the case of 4D was considered, but that it is easily generalized for any dimensions. 
 
 Any oscillatory solution can be obtained as the circle Fourier transform \cite{Wolf89} (page 115):
 
 \begin{proposition}
  Any oscillatory solution $H\in {\cal F}_\lambda$ of eq. (\ref{EqHelmholtz}) can be written as:
  \begin{equation}
   H(\varphi,v,w) =\int_{-\pi}^\pi d\phi e^{i\lambda(v\sin\phi + w \cos\phi)}h(\phi,\varphi)
   \label{circleFourier}
  \end{equation}
  \end{proposition}
  
  \noindent \textit{Proof:} See \cite{Wolf89}, page 115.
  
  In particular, functions $F(\varphi,v,w)=\langle \eta_{(\varphi,v,w)},f\rangle$ are of this type, with $h(\phi,\varphi)=\eta(\phi-\varphi)\overline{f(\phi)}$ and also $\eta_{(\varphi,v,w)}(\theta)$, with $h(\phi,\varphi)=\eta(\phi-\varphi)\delta(\phi-\theta)$.

 As a consequence of Prop.  \ref{NonAdmissible}, even if $h(\cdot,\varphi)\in{\cal H}$, $H(\varphi,v,w)$ is not square integrable on $E(2)$\footnote{Note that the circle Fourier transform (\ref{circleFourier}) is not an ordinary 2D Fourier transform, and therefore standard results on Fourier analysis do not apply here. This implies that the circle Fourier transform is not unitary with respect to the left invariant Haar measure on $E(2)$.}.
 
 To better understand the root of the problem, it is convenient to study the properties of oscillatory solutions of eq. (\ref{EqHelmholtz}).
 
 For oscillatory solutions, the Cauchy initial value problem for Helmholtz equation is well defined (see \cite{Wolf81,SphereMomentum20}).

\begin{proposition}\label{PropIVP}
   
  The initial value problem for Helmholtz equation:

\begin{equation}\label{IVPHemlholtz}
 \left(\frac{\partial^2\,}{\partial v^2}+\frac{\partial^2\,}{\partial w^2}+\lambda^2\right)F(v,w)=0 \,,\quad F(v,w)|_{w=0}=\Fc(v)\,,\qquad \left.\frac{\partial F(v,w)}{\partial w}\right|_{w=0}=\Fb(v)
\end{equation}
is well-posed, provided that the Fourier transform of the functions $\Fc(v)$ and  $\Fb(v)$, the initial data for the function  $F(v,w)$ and the normal derivative in the line $w=0$, respectively, have support on the interval $[-\lambda,\lambda]$.
\end{proposition}

\begin{remark}\label{RemarkInvariance}
 Prop. \ref{PropIVP}  is invariant under the Euclidean group, in the sense that the initial condition can be taken in any line of the plane $(v,w)$.
\end{remark}

The solution $F(v,w)$ can be \textit{evolved} from the initial data by convolution with certain integral kernels (see Appendix B of \cite{SphereMomentum20}).
\begin{proposition}\label{PropPropagationHelmholtz}
Any oscillatory solution $F(v,w)$ of Eq. (\ref{IVPHemlholtz})  can be obtained from the initial 
data $\Fc(v)$ and  $\Fb(v)$ by:

 \begin{eqnarray}
 F(v,w)&=&\int_\mathbb{R} dv'\left[ \Delta(v-v',w-w')\frac{\partial F(v',w')}{\partial w'}-\frac{\partial \Delta(v-v',w-w')}{\partial w'}F(v',w')\right]_{w'=0} \nonumber \\
  &=& \int_\mathbb{R} dv'\left(\Delta(v-v',w) \Fb(v') + \frac{\partial \Delta(v-v',w)}{\partial w} \Fc(v')\right) \label{PropagationHelmholtz}
\end{eqnarray}
where 
\begin{equation}
 \Delta(v,w)=\frac{1}{\sqrt{2\pi}}\int_{-\lambda}^\lambda d\epsilon  e^{i \epsilon v} \frac{\sin(\sqrt{\lambda^2-\epsilon^2}w)}{\sqrt{\lambda^2-\epsilon^2}} = \frac{1}{\sqrt{2\pi}}\int_{-\pi/2}^{\pi/2}d\theta
 e^{i \lambda v \sin\theta}\sin(\lambda w \cos\theta) 
\end{equation}
%
%
with $\Delta(v,w)$ a solution of the IVP  (\ref{IVPHemlholtz}) with initial data:
\begin{eqnarray}
 \Dc(v) &=& \Delta(v,0)=0\\
 \Db(v)& =& \frac{\partial \Delta(v,w)}{\partial w}|_{w=0} = \frac{1}{\sqrt{2\pi}}\int_{-\lambda}^\lambda d\epsilon  e^{i \epsilon v} = \sqrt{\frac{2}{\pi}} \frac{\sin(\lambda v)}{v} \equiv \sqrt{\frac{2}{\pi}}\lambda k_{1/2}( \lambda v)
\end{eqnarray}
where $k_\nu(x)\equiv \frac{J_\nu(x)}{x^\nu}$ are entire functions known as Bochner-Riesz integral convolution kernels.

\end{proposition}

\noindent\text{Proof:} See Appendix B of \cite{SphereMomentum20} (translated to the case of $E(2)$) and apply convolution theorem to eq. (B.10) of \cite{SphereMomentum20}.

Using Prop. \ref{PropPropagationHelmholtz}, it can be easily proven that for any oscillatory solution $F(\varphi,v,w)$, the $L^2$-norm  with respect to $w$ diverges:
\begin{equation}
\int_{-\infty}^\infty dw |F(\varphi,v,w)|^2 =+\infty
\end{equation}

Remark \ref{RemarkInvariance} also applies here, therefore the $L^2$-norm of $F(\varphi,v,w)$
along any line in the $(v,w)$ plane diverges. The reason of this divergence is traced back to the oscillatory behaviour of $F(\varphi,v,w)$ with respect to $(v,w)$. And this divergence is responsible for the divergence in Proposition \ref{NonAdmissible}.

According to this, it is clear that any resolution of the identity with the Haar measure involving translations will diverge, even if we restrict to one of the translations variable, unless an admissibility condition is imposed to $\eta$, as is the case  of de Bi\`evre's proposal \cite{DeBievre89}  (see also  \cite{Torresani95}).

We will introduce a new resolution of the identity for coherent states of the $E(2)$ that involves an integration suited to solutions of Helmholtz equation, that restricts the integration to the initial data of Proposition \ref{PropIVP}. But first we will discuss the already known inner product for oscillatory solutions of Helmholtz equation \cite{Wolf81}.

\section{An invariant inner product for solutions of the Helmholtz equation}

We know  a scalar product invariant under $\E(2)$ for solutions of the
Helmholtz equation is given by \cite{Wolf81}:
\begin{eqnarray}
  \langle f,g\rangle
  &=\int_{\R^2}dv\,dv'\,\Big[
    K_1(v-v')\,\conj{f(v,w)}\,g(v',w)\nonumber\\
  &\qquad+K_2(v-v')\,\conj{\partial_w f(v,w)}\,\partial_w g(v',w)\Big]
\end{eqnarray}
with  $K_1(x)=\dfrac{J_1(x)}{x}$, \quad $K_2(x)=J_0(x)$.

This scalar product is invariant under the group $\E(2)$ \cite{Wolf81}, in particular
it is invariant under $T_w$, and therefore it does not depend on the value of $w$, thus we choose without lack of generality $w=0$:
\begin{eqnarray}\label{innerHelmholtz}
  \langle f,g\rangle &=&\int_{\R^2}dv\,dv'\,\Big[
    K_1(v-v')\,\conj{f(v,0)}\,g(v',0)\nonumber\\
  &\quad&+K_2(v-v')\,\partial_w \conj{f(v,0)}\,\partial_w g(v',0)\Big]
\end{eqnarray}
where $f(v,0)$ and $\partial_w f(v,0)$ are the initial data
for $f$ and its derivative $\partial_w f$ at $w=0$ that uniquely determine a
solution of the Helmholtz equation (see Proposition \ref{PropIVP}).

\section{A new resolution of the identity on $E(2)$ adapted to solutions of Helmholtz equation}

Using Proposition \ref{Aidentity} and Remark \ref{Remark-ResolutionIdentity}, let us propose the following resolution of the identity for coherent states of the $E(2)$, using the inner product for solutions of the Helmholtz equation given by eq. (\ref{innerHelmholtz}).

\begin{theorem}
 Suppose that $F(\varphi,v,w)=\langle \eta_{(\varphi,v,w)},f\rangle$ and $G(\varphi,v,w)=\langle \eta_{(\varphi,v,w)},g\rangle$. Since $F$ and $G$ satisfy Helmholtz equation,  we have the following resolution of the identity:
 \begin{eqnarray}
\langle F,G\rangle &=& \int_{S^1}d\varphi\int_{\R^2}dv\,dv'\,\Big[K_1(v-v')\,\conj{F(\varphi,v,0)}\,G(\varphi,v',0)\nonumber\\
  & &\qquad\qquad+K_2(v-v')\,\partial_w \conj{F(\varphi,v,0)}\,\partial_w G(\varphi,v',0)\Big] = \langle f, \tilde{A}_{\eta}g\rangle = \langle f, g\rangle
\end{eqnarray}

\end{theorem}

\noindent\textit{Proof:} The proof is given in Appendix \ref{apResId}, which uses results proven in previous appendices.

\begin{remark}
 Note that the operator $\tilde{A}_{\eta}$ is similar to the operator $A_{\eta}$ introduced in Proposition \ref{Abounded}, but instead of using the Haar measure in $E(2)$, the integration is changed by that of eq. (\ref{innerHelmholtz}) (since the involved functions  satisfy Helmholtz equation), adding the integration on rotations with its corresponding Haar measure. This implies that the adjoint $T_\eta^{\dag}$ given in Proposition \ref{adjointCSmap} must be replaced by $\tilde{T}_\eta^{\dag}$, using the same integration as $\tilde{A}_{\eta}$.
 
\end{remark}

\section{Resolution of the Identity without integrating over $\varphi$}

Note that if we do not integrate over $\varphi$, and consider only the states
of the form $\eta_{(0,v,w)}(\theta)$,
\[
  A_{m,n}=\frac{1}{2\pi}\int_{-\pi}^{\pi}d\theta\;e^{i(n-m)\theta}\,|\eta(\theta)|^2.
\]
In this case it is enough to take $\eta(\theta)=1$ to obtain
\[
  A_{m,n}=\delta_{m,n}\quad\Rightarrow\quad A=I.
\]

\section{Recovering the result of De Bi\`evre and Torr\'esani}

We can also recover the result of De Bi\`evre and Torr\'esani by setting
$K_1(x)=\delta(x)$, $K_2(x)=0$, taking $w=0$ and using
$\tilde\eta_{(\varphi,v,0)}(\theta)$.

If we use, in addition,
\begin{align*}
  \langle\tilde\eta_{(\varphi,v,0)},f\rangle=\tilde f(\varphi,v,0)
  &=\int_{-\pi}^{\pi}d\theta\;e^{-iv\sin(\theta-\varphi)}\,
    \conj{\eta(\theta-\varphi)}\,f(\theta)\\
  &=\int_{-\pi}^{\pi}d\theta\;e^{-iv\sin\theta}\,\conj{\eta(\theta)}\,f(\theta+\varphi)
\end{align*}

In this case $\hat K_1(\sin\theta')=1$, \ $\hat K_2(\sin\theta)=0$.

We have
\begin{align*}
  A_{m,n}&=\frac{1}{2\pi}\int_{S^1}d\varphi\;e^{i(n-m)\varphi}
    \int_{-\pi}^{\pi}d\theta\;e^{i(n-m)\theta}
    \left(\frac{\hat K_1(\sin\theta)}{|\cos\theta|}+|\cos\theta|\,\hat K_2(\sin\theta)\right)
    |\eta(\theta)|^2\\
  &\quad+\frac{1}{2\pi}\,e^{-im\pi}\int_{S^1}d\varphi\;e^{i(n-m)\varphi}
    \int_{-\pi}^{\pi}d\theta\;e^{i(n+m)\theta}
    \left(\frac{\hat K_1(\sin\theta)}{|\cos\theta|}-|\cos\theta|\,\hat K_2(\sin\theta)\right)
    \conj{\eta(\theta)}\,\eta(\pi-\theta)
\end{align*}
\[
  \Rightarrow\quad A_{m,n}=\delta_{m,n}\int_{-\pi}^{\pi}d\theta\;
  \frac{|\eta(\theta)|^2}{|\cos\theta|}
  +e^{-im\pi}\,\delta_{m,n}\int_{-\pi}^{\pi}d\theta\;e^{2in\theta}\,
  \frac{\conj{\eta(\theta)}\,\eta(\pi-\theta)}{|\cos\theta|}
\]
If we want $A_{m,n}\propto\delta_{m,n}$ we must take $\eta(\theta)$ with
support in $(-\frac{\pi}{2},\frac{\pi}{2})$ and satisfying
\[
  \int_{-\pi/2}^{\pi/2}d\theta\;\frac{|\eta(\theta)|^2}{\cos\theta}<\infty
\]
and in this way we recover the result of De Bi\`evre and Torr\'esani.

\section{Conclusions}

In this paper we have constructed a non-divergent resolution of the Identity for Perelomov coherent states of $E(2)$,  based in a non-local integration measure designed for solutions of Helmholtz equation in \cite{Wolf81}. Making a different choice of the integration measure we have recovered the admissibility condition given in \cite{DeBievre89}.

The construction of a non-divergent resolution of the Identity for the case of $E(n)$ is stright-forward using the same technique. This work is under progress and will be published elsewhere. Also, this construction could be applied to the case of the Poincar\'e group, using an invariant non-local integration with suitable convolution kernels, adapted to solutions of the Klein-Gordon equation (see \cite{LopezRuiz20,LopezRuiz21}). This will also be published elsewhere.

%
%
%

%

\section*{Acknowledgements}

The author thanks the support of the Spanish MICINN through the project PID2022-138144NB-I00. I wish to give special thanks to G. Garrig\'os for enlightling discussions.

\appendix

\section{Proof of Proposition \ref{NonAdmissible}}
\label{app:0}

Let us compute a representation of the resolution operator $A_\eta$ in a suitable basis, such as for
example
\[
  \{f_n(\theta)\}_{n\in\Z},\qquad f_n(\theta)=\frac{1}{\sqrt{2\pi}}\,e^{in\theta}.
\]
\begin{align*}
  (A_\eta)_{m,n}=\langle f_m,Af_n\rangle
  &=\frac{1}{2\pi}\int_{-\pi}^{\pi}d\theta'\,e^{-im\theta'}
    \int_{\E(2)}d\mu(\varphi,v,w)\int_{-\pi}^{\pi}d\theta\;
    e^{-i(v\sin\theta+w\cos\theta)}\,\conj{\eta(\theta-\varphi)}\,e^{in\theta}\\
  &\qquad\cdot e^{i(v\sin\theta'+w\cos\theta')}\,\eta(\theta'-\varphi)\\
  &=\frac{1}{2\pi}\int_{-\pi}^{\pi}d\theta'\,e^{-im\theta'}
    \int_{-\pi}^{\pi}d\theta\,e^{in\theta}
    \int_{\E(2)}d\mu(\varphi,v,w)\;
    e^{i[v(\sin\theta'-\sin\theta)+w(\cos\theta'-\cos\theta)]}\\
  &\qquad\cdot\conj{\eta(\theta-\varphi)}\,\eta(\theta'-\varphi).
\end{align*}
If we take the Haar measure on $\E(2)$,
$d\mu(\varphi,v,w)=d\varphi\,dv\,dw$,
\begin{align*}
  \Rightarrow\quad A_{m,n}
  &=\frac{1}{2\pi}\int_{-\pi}^{\pi}d\theta'\,e^{-im\theta'}
    \int_{-\pi}^{\pi}d\theta\,e^{in\theta}
    \int_{\R^2}dv\,dw\;
    e^{i[v(\sin\theta'-\sin\theta)+w(\cos\theta'-\cos\theta)]}\\
  &\qquad\cdot\int_{S^1}d\varphi\;\conj{\eta(\theta-\varphi)}\,\eta(\theta'-\varphi).
\end{align*}
Let us compute
\begin{align*}
  B(\theta,\theta')
  &=\int_{\R^2}dv\,dw\;e^{i[v(\sin\theta'-\sin\theta)+w(\cos\theta'-\cos\theta)]}\\
  &=\int_{\R}e^{iv(\sin\theta'-\sin\theta)}\,dv
    \int_{\R}e^{iw(\cos\theta'-\cos\theta)}\,dw.
\end{align*}
Using  Appendix \ref{app:2}:
\begin{align*}
  B(\theta,\theta')
  &=\left(\frac{1}{|\cos\theta|}\,\delta(\theta'-\theta)
         +\frac{1}{|\cos\theta|}\,\delta(\theta'-(\pi-\theta))\right)\\
  &\qquad\cdot\left(\frac{1}{|\sin\theta|}\,\delta(\theta'-\theta)
         +\frac{1}{|\sin\theta|}\,\delta(\theta'+\theta)\right)\\
  &=\frac{1}{\tfrac12|\sin2\theta|}\Big[\delta(\theta'-\theta)^2
     +\delta(\theta'-(\pi-\theta))\,\delta(\theta'+\theta)\\
  &\qquad+\delta(\theta'-\theta)\,\delta(\theta'+\theta)
     +\delta(\theta'-(\pi-\theta))\,\delta(\theta'-\theta)\Big]
\end{align*}

\begin{align*}
  &\theta'+\theta=\pi\ \text{ and }\ \theta'+\theta=0 &&\text{Impossible}\\
  &\theta'=\theta\ \text{ and }\ \theta'+\theta=0 &&\Rightarrow\ \theta=\theta'=0\\
  &\theta'+\theta=\pi\ \text{ and }\ \theta'=\theta &&\Rightarrow\ \theta=\theta'=\tfrac{\pi}{2}
\end{align*}
Note that in any case we are going to obtain an infinite result for
$(A_\eta)_{m,n}$, since $\delta(\theta'-\theta)^2$ will give infinity, independently of
the chosen $\eta(\theta)$.

This means that we cannot use the Haar measure and integrate over the whole
group.

\section{Dirac delta of a function with simple zeros}
\label{app:1}

We show that if $f\colon\R\to\R$ has $n$ simple zeros $x_1,\dots,x_n$,
\[
  f(x_i)=0,\qquad i=1,\dots,n,
\]
and is $C^1$ on $\R$, then
\begin{equation}\label{eq:deltaf}
  \delta\bigl(f(x)\bigr)=\sum_{i=1}^{n}\frac{1}{|f'(x_i)|}\,\delta(x-x_i).
\end{equation}

To prove it, we first consider the case $n=1$:
\[
  \delta\bigl(f(x)\bigr)\quad\text{with}\quad f(x_1)=0,\qquad f'(x_1)\neq0.
\]
If we approximate $\delta(x)$ by
\[
  \eta_\varepsilon(x)=\frac{1}{\varepsilon\sqrt{\pi}}\,e^{-x^2/\varepsilon^2},
  \qquad \varepsilon>0,
  \qquad
  \lim_{\varepsilon\to0}\eta_\varepsilon(x)=\delta(x),
\]
then
\[
  \eta_\varepsilon\bigl(f(x)\bigr)=\frac{1}{\varepsilon\sqrt{\pi}}\,
  e^{-f(x)^2/\varepsilon^2}.
\]
Clearly, $\eta_\varepsilon\bigl(f(x)\bigr)$ is significantly different from
zero only when $x\simeq x_1$, so in this region we may approximate
$f(x)\simeq f'(x_1)(x-x_1)$:
\[
  \eta_\varepsilon\bigl(f(x)\bigr)\simeq\frac{1}{\varepsilon\sqrt{\pi}}\,
  e^{-f'(x_1)^2(x-x_1)^2/\varepsilon^2}
  =\frac{1}{|f'(x_1)|}\,\eta_{\varepsilon'}(x-x_1),
  \qquad\text{with}\quad \varepsilon'=\frac{\varepsilon}{|f'(x_1)|}>0.
\]
In the limit $\varepsilon\to0$ we have $\varepsilon'\to0$ and
$\eta_{\varepsilon'}(x-x_1)\to\delta(x-x_1)$, hence
\[
  \eta_\varepsilon\bigl(f(x)\bigr)\longrightarrow\frac{1}{|f'(x_1)|}\,\delta(x-x_1).
\]

When $f(x)$ has more than one zero, a similar argument shows that
\[
  \eta_\varepsilon\bigl(f(x)\bigr)\simeq\frac{1}{\varepsilon\sqrt{\pi}}
  \sum_{i=1}^{n}e^{-f'(x_i)^2(x-x_i)^2/\varepsilon^2}
  =\sum_{i=1}^{n}\frac{1}{|f'(x_i)|}\,\eta_{\varepsilon_i'}(x-x_i),
  \qquad\text{with}\quad \varepsilon_i'=\frac{\varepsilon}{|f'(x_i)|},
  \quad i=1,\dots,n.
\]
Letting $\varepsilon\to0$, we obtain
\[
  \delta\bigl(f(x)\bigr)=\sum_{i=1}^{n}\frac{1}{|f'(x_i)|}\,\delta(x-x_i).
\]

\subsection{Generalization to two variables: $h(x,y)=f(x)\,g(y)$}
\label{app:1-2d}

\noindent\textit{(This subsection is new: it is not in the manuscript.)}

Let $f,g\colon\R\to\R$ be $C^1$, where the only zeros of $f$ are
$x_1,\dots,x_n$ and the only zeros of $g$ are $y_1,\dots,y_m$, all of them
simple:
\[
  f(x_i)=0,\quad f'(x_i)\neq0,\qquad g(y_j)=0,\quad g'(y_j)\neq0 .
\]
The zero set of $h(x,y)=f(x)g(y)$ in $\R^2$ is the union of the vertical lines
$x=x_i$ and the horizontal lines $y=y_j$. These lines cross at the $nm$ points
$P_{ij}=(x_i,y_j)$, which are exactly the points where
$\nabla h=\bigl(f'(x)g(y),\,f(x)g'(y)\bigr)$ vanishes.

\begin{proposition}
For every $\varphi\in C_c(\R^2)$ whose support contains no crossing point
$P_{ij}$,
\[
  \lim_{\varepsilon\to0}\int_{\R^2}\eta_\varepsilon\bigl(f(x)g(y)\bigr)\,
  \varphi(x,y)\,dx\,dy
  =\sum_{i=1}^{n}\int_\R\frac{\varphi(x_i,y)}{|f'(x_i)|\,|g(y)|}\,dy
  +\sum_{j=1}^{m}\int_\R\frac{\varphi(x,y_j)}{|f(x)|\,|g'(y_j)|}\,dx ,
\]
that is, on $\R^2\setminus\{P_{ij}\}$,
\begin{equation}\label{eq:delta2d}
  \delta\bigl(f(x)g(y)\bigr)
  =\sum_{i=1}^{n}\frac{1}{|f'(x_i)\,g(y)|}\,\delta(x-x_i)
  +\sum_{j=1}^{m}\frac{1}{|f(x)\,g'(y_j)|}\,\delta(y-y_j).
\end{equation}
\end{proposition}

\noindent\textit{Proof:}
\emph{Step 1: separate the two families of lines.} Let
$U=\{(x,y):g(y)\neq0\}$ and $V=\{(x,y):f(x)\neq0\}$. Both are open and
$U\cup V=\R^2\setminus\{P_{ij}\}$. Since the support of $\varphi$ is compact
and contained in $U\cup V$, a partition of unity gives
$\varphi=\varphi_1+\varphi_2$, with $\operatorname{supp}\varphi_1\subset U$ and
$\operatorname{supp}\varphi_2\subset V$, both compact. In particular there is
a $c>0$ with $|g(y)|\ge c$ on $\operatorname{supp}\varphi_1$ and
$|f(x)|\ge c$ on $\operatorname{supp}\varphi_2$.

\emph{Step 2: the term with $\varphi_1$.} For each fixed $y$ with
$(x,y)\in\operatorname{supp}\varphi_1$, the function
$x\mapsto F_y(x)=f(x)g(y)$ has simple zeros $x_1,\dots,x_n$, with
$F_y'(x_i)=f'(x_i)g(y)$. As in Appendix~\ref{app:1}, near $x_i$ we have
$F_y(x)\simeq f'(x_i)g(y)(x-x_i)$ and
\[
  \eta_\varepsilon\bigl(f(x)g(y)\bigr)\simeq
  \frac{1}{|f'(x_i)\,g(y)|}\,\eta_{\varepsilon_i(y)}(x-x_i),
  \qquad
  \varepsilon_i(y)=\frac{\varepsilon}{|f'(x_i)\,g(y)|}\le\frac{\varepsilon}{c\,|f'(x_i)|}.
\]
The bound $|g|\ge c$ makes the width $\varepsilon_i(y)$ go to zero uniformly
in $y$. Hence the $x$-integral is bounded uniformly in $y$ and converges to
$\sum_i\varphi_1(x_i,y)/|f'(x_i)g(y)|$. By dominated convergence,
\[
  \int\eta_\varepsilon(h)\,\varphi_1\;\longrightarrow\;
  \sum_{i=1}^{n}\int_\R\frac{\varphi_1(x_i,y)}{|f'(x_i)|\,|g(y)|}\,dy .
\]
The lines $y=y_j$ contribute nothing here, because they do not meet
$\operatorname{supp}\varphi_1$.

\emph{Step 3: the term with $\varphi_2$.} This is the same argument with the
roles of $x$ and $y$ exchanged. It gives
$\sum_j\int\varphi_2(x,y_j)/|f(x)g'(y_j)|\,dx$, and the lines $x=x_i$
contribute nothing.

\emph{Step 4.} It remains to add the two terms. Since $\varphi_2(x_i,\cdot)=0$
on the lines $x=x_i$ and $\varphi_1(\cdot,y_j)=0$ on the lines $y=y_j$, we have
$\varphi_1(x_i,y)=\varphi(x_i,y)$ and $\varphi_2(x,y_j)=\varphi(x,y_j)$. This
gives the result.$\square$

The factor $1/|f'(x_i)g(y)|$ is $1/|\nabla h|$ on the line $x=x_i$, and
$1/|f(x)g'(y_j)|$ is $1/|\nabla h|$ on the line $y=y_j$. Hence
\eqref{eq:delta2d} is a special case of the general formula
$\delta(h)=\delta_{\{h=0\}}/|\nabla h|$ (surface measure on the zero set
divided by the gradient).

\paragraph{Crossing points.} The assumption that the support of $\varphi$
avoids the points $P_{ij}$ is necessary. Near $P_{ij}$, with $s=x-x_i$,
$t=y-y_j$, $a=f'(x_i)$, $b=g'(y_j)$, we have $h\simeq ab\,st$. If
$\varphi\equiv1$ on the square $|s|,|t|\le L$, then
\[
  \int_{-L}^{L}\!\!\int_{-L}^{L}\eta_\varepsilon(ab\,st)\,ds\,dt
  =\frac{2}{|ab|}\int_0^{|ab|L^2/\varepsilon}\frac{\operatorname{erf}\tau}{\tau}\,d\tau
  =\frac{2}{|ab|}\,\ln\frac{1}{\varepsilon}+O(1),
\]
which diverges logarithmically. The same can be seen in \eqref{eq:delta2d}:
the density $1/|g(y)|\sim1/|b(y-y_j)|$ is not integrable at $y_j$. Therefore
$\delta(f(x)g(y))$ does not define a distribution on test functions that do
not vanish at the $P_{ij}$. This is the same phenomenon that appears in
Appendix~\ref{app:2} at $x=\pi/2$, where $u=v=0$ simultaneously and the
coefficient $1/|\cos x|$ diverges.

\paragraph{Relation to Appendix~\ref{app:2}.} Formula \eqref{eq:deltafg} is
\eqref{eq:delta2d} with $(u,v)$ in place of $(x,y)$ and a single zero in each
variable. Appendix~\ref{app:2} uses it along the line
$x'\mapsto(u,v)=(x-x',\,\pi-x-x')$. This line crosses each line of zeros
transversally, with $|\partial u/\partial x'|=|\partial v/\partial x'|=1$,
which is why the coefficients agree with those of the one-dimensional formula
in $x'$.

\section{Delta of $\sin\theta-\sin\theta'$ and of $\cos\theta-\cos\theta'$}
\label{app:2}

We now consider $\delta\bigl(B(x,x')\bigr)$, specifically the case
\[
  \delta\bigl(f(x-x')\,g(x+x')\bigr)\qquad\text{with}\quad f(0)=g(0)=0.
\]
Let
\[
  u=x-x',\qquad v=\pi-(x+x')
  \quad\Longrightarrow\quad
  \delta\bigl(f(u)\,g(v)\bigr).
\]

Applying Appendix~\ref{app:1}:
\begin{equation}\label{eq:deltafg}
  \delta\bigl(f(u)\,g(v)\bigr)=\frac{1}{|f'(0)\,g(v)|}\,\delta(u)
  +\frac{1}{|f(u)\,g'(0)|}\,\delta(v).
\end{equation}


\subsection{Case $\sin x-\sin x'$}

Now take
\[
  \sin x-\sin x'=2\sin\Bigl(\frac{x-x'}{2}\Bigr)\cos\Bigl(\frac{x+x'}{2}\Bigr),
\]
\[
  f(u)=2\sin\frac{u}{2}=2\sin\Bigl(\frac{x-x'}{2}\Bigr),\qquad
  g(v)=\sin\frac{v}{2}=\sin\Bigl(\frac{\pi}{2}-\frac{x+x'}{2}\Bigr)
  =\cos\frac{x+x'}{2}.
\]
Then
\begin{align*}
  \delta\Bigl(2\sin\frac{u}{2}\,\sin\frac{v}{2}\Bigr)
  &=\frac{1}{\bigl|\cos0\cdot\sin\frac{v}{2}\bigr|}\,\delta(u)
   +\frac{1}{\bigl|\sin\frac{u}{2}\cdot\cos0\bigr|}\,\delta(v)\\
  &=\frac{1}{\bigl|\sin\frac{v}{2}\bigr|}\,\delta(u)
   +\frac{1}{\bigl|\sin\frac{u}{2}\bigr|}\,\delta(v).
\end{align*}
Since
\begin{align*}
  u=0&\;\Rightarrow\; x=x' \;\Rightarrow\; v=\pi-2x,\\
  v=0&\;\Rightarrow\; x+x'=\pi \;\Rightarrow\; x'=\pi-x
     \;\Rightarrow\; u=x-(\pi-x)=2x-\pi,
\end{align*}
we get
\begin{align*}
  \delta\Bigl(2\sin\frac{u}{2}\,\sin\frac{v}{2}\Bigr)
  &=\frac{1}{\bigl|\sin\bigl(\frac{\pi}{2}-x\bigr)\bigr|}\,\delta(x'-x)
   +\frac{1}{\bigl|\sin\bigl(x-\frac{\pi}{2}\bigr)\bigr|}\,\delta\bigl(x'-(\pi-x)\bigr)\\
  &=\frac{1}{|\cos x|}\,\delta(x'-x)+\frac{1}{|\cos x|}\,\delta\bigl(x'-(\pi-x)\bigr).
\end{align*}
Therefore,
\begin{equation}\label{eq:deltasen}
  \delta(\sin\theta-\sin\theta')=\frac{1}{|\cos\theta|}\,\delta(\theta-\theta')
  +\frac{1}{|\cos\theta|}\,\delta\bigl(\theta'-(\pi-\theta)\bigr).
\end{equation}

Note that according to the comment on crossingpoints in Appendix \ref{app:1}, $\delta(\sin\theta-\sin\theta')$ is a true distribution on functions whose support does not include 
$\theta=\frac{\pi}{2}, \frac{3\pi}{2}$.

\subsection{Case $\cos x-\cos x'$}

If we take
\[
  \cos x-\cos x'=-2\sin\Bigl(\frac{x-x'}{2}\Bigr)\sin\Bigl(\frac{x+x'}{2}\Bigr),
\]
\[
  f(u)=-2\sin\frac{u}{2}=-2\sin\Bigl(\frac{x-x'}{2}\Bigr),\qquad
  g(v)=\sin\frac{v}{2}=\sin\Bigl(\frac{x+x'}{2}\Bigr),
\]
\begin{align*}
  \delta\Bigl(-2\sin\frac{u}{2}\,\sin\frac{v}{2}\Bigr)
  &=\frac{1}{\bigl|-\cos0\cdot\sin\frac{v}{2}\bigr|}\,\delta(u)
   +\frac{1}{\bigl|-\sin\frac{u}{2}\cdot\cos0\bigr|}\,\delta(v)\\
  &=\frac{1}{\bigl|\sin\frac{v}{2}\bigr|}\,\delta(u)
   +\frac{1}{\bigl|\sin\frac{u}{2}\bigr|}\,\delta(v).
\end{align*}
Since
\begin{align*}
  u=0&\;\Rightarrow\; x=x' \;\Rightarrow\; v=2x,\\
  v=0&\;\Rightarrow\; x+x'=0 \;\Rightarrow\; x'=-x
     \;\Rightarrow\; u=x-(-x)=2x,
\end{align*}
we get
\[
  \delta\Bigl(-2\sin\frac{u}{2}\,\sin\frac{v}{2}\Bigr)
  =\frac{1}{|\sin x|}\,\delta(x-x')+\frac{1}{|\sin x|}\,\delta(x+x').
\]
Therefore,
\begin{equation}\label{eq:deltacos}
  \delta(\cos\theta-\cos\theta')=\frac{1}{|\sin\theta|}\,\delta(\theta-\theta')
  +\frac{1}{|\sin\theta|}\,\delta(\theta+\theta').
\end{equation}

Note that according to the comment on crossing points in Appendix \ref{app:1}, $\delta(\cos\theta-\cos\theta')$ is a true distribution on functions whose support does not include 
$\theta=0, \pi$.

\section{Computation of a double integral}
\label{app:3}

We compute
\[
  \int_{\R\times\R}dv\,dv'\,k(v-v')\,e^{-i(v\sin\theta-v'\sin\theta')}
  =\int_{\R\times\R}dv\,dw\,k(w)\,e^{-i(v\sin\theta-(v-w)\sin\theta')},
\]
where we set
\[
  \begin{cases} v=v\\ w=v-v'\end{cases}
  \qquad
  \begin{cases} v=v\\ v'=v-w\end{cases}
  \qquad
  \frac{\partial(v,v')}{\partial(v,w)}=
  \begin{vmatrix}1&0\\1&-1\end{vmatrix}=-1.
\]
Then
\begin{align*}
  &=\int_{\R}dw\,k(w)\,e^{-iw\sin\theta'}\int_{\R}dv\,e^{-iv(\sin\theta-\sin\theta')}
  \qquad\text{(using Appendix~\ref{app:2})}\\
  &=2\pi\hat k(\sin\theta')\,\frac{1}{|\cos\theta|}
    \Bigl(\delta(\theta-\theta')+\delta\bigl(\theta'-(\pi-\theta)\bigr)\Bigr)\,,
\end{align*}
Since $\sin(\pi-\theta)=\sin\theta$,
\[
  =2\pi\frac{\hat k(\sin\theta)}{|\cos\theta|}
  \Bigl(\delta(\theta-\theta')+\delta\bigl(\theta'-(\pi-\theta)\bigr)\Bigr).
\]

Note that this is a distribution on the whole interval $(-\pi,\pi)$ if $\hat{k}(\pm 1)=0$.

\section{Computation of $\tilde{A}_\eta$ with the invariant scalar product}\label{apResId}

Define $f_m(\theta)=\frac{1}{\sqrt{2\pi}}e^{im\theta}$. Then

\begin{align*}
  (\tilde{A}_{\eta})_{m,n}&=\frac{1}{2\pi}\langle f_m,Af_n\rangle\\
  &=\frac{1}{2\pi}\int_{-\pi}^{\pi}d\theta'\,e^{-im\theta'}
    \int_{S^1}d\varphi\int_{\R\times\R}dv\,dv'\,\Big[\\
  &\qquad K_1(v-v')\int_{-\pi}^{\pi}d\theta\;e^{-iv\sin\theta}\,
    \conj{\eta(\theta-\varphi)}\,e^{in\theta}\,e^{iv'\sin\theta'}\,\eta(\theta'-\varphi)\\
  &\qquad+K_2(v-v')\int_{-\pi}^{\pi}d\theta\;\cos\theta\;e^{-iv\sin\theta}\,
    \conj{\eta(\theta-\varphi)}\,e^{in\theta}\,\cos\theta'\,e^{iv'\sin\theta'}\,
    \eta(\theta'-\varphi)\Big]\\
  &=\frac{1}{2\pi}\int_{-\pi}^{\pi}d\theta'\,e^{-im\theta'}
    \int_{-\pi}^{\pi}d\theta\,e^{in\theta}
    \int_{S^1}d\varphi\;\conj{\eta(\theta-\varphi)}\,\eta(\theta'-\varphi)\\
  &\qquad\cdot\int_{\R\times\R}dv\,dv'\,\Big[
    K_1(v-v')\,e^{-i(v\sin\theta-v'\sin\theta')}
    +\cos\theta\cos\theta'\,K_2(v-v')\,e^{-i(v\sin\theta-v'\sin\theta')}\Big]
\end{align*}
\begin{center}(see previous appendices)\end{center}
\begin{align*}
  &=\frac{1}{2\pi}\int_{-\pi}^{\pi}d\theta'\,e^{-im\theta'}
    \int_{-\pi}^{\pi}d\theta\,e^{in\theta}
    \int_{S^1}d\varphi\;\conj{\eta(\theta-\varphi)}\,\eta(\theta'-\varphi)\\
  &\qquad\cdot\Big[\hat K_1(\sin\theta')\Big(\frac{1}{|\cos\theta|}\,\delta(\theta-\theta')
    +\frac{1}{|\cos\theta|}\,\delta(\theta'-(\pi-\theta))\Big)\\
  &\qquad\quad+\cos\theta\cos\theta'\,\hat K_2(\sin\theta')
    \Big(\frac{1}{|\cos\theta|}\,\delta(\theta-\theta')
    +\frac{1}{|\cos\theta|}\,\delta(\theta'-(\pi-\theta))\Big)\Big]\\[4pt]
  &=\frac{1}{2\pi}\int_{-\pi}^{\pi}d\theta\;e^{i(n-m)\theta}\,
    \frac{\hat K_1(\sin\theta)}{|\cos\theta|}\int_{S^1}d\varphi\;|\eta(\theta-\varphi)|^2\\
  &\quad+\frac{1}{2\pi}\int_{-\pi}^{\pi}e^{-im\pi}\,e^{i(n+m)\theta}\,
    \frac{\hat K_1(\sin\theta)}{|\cos\theta|}\int_{S^1}d\varphi\;
    \conj{\eta(\theta-\varphi)}\,\eta(\pi-\theta-\varphi)\\
  &\quad+\frac{1}{2\pi}\int_{-\pi}^{\pi}e^{i(n-m)\theta}\,|\cos\theta|\,
    \hat K_2(\sin\theta)\int_{S^1}d\varphi\;|\eta(\theta-\varphi)|^2\\
  &\quad-\frac{1}{2\pi}\int_{-\pi}^{\pi}d\theta\;e^{-im\pi}\,e^{i(n+m)\theta}\,
    |\cos\theta|\,\hat K_2(\sin\theta)\cdot\int_{S^1}d\varphi\;
    \conj{\eta(\theta-\varphi)}\,\eta(\pi-\theta-\varphi)\\[4pt]
  &=\frac{1}{2\pi}\int_{-\pi}^{\pi}d\theta\;e^{i(n-m)\theta}
    \left(\frac{\hat K_1(\sin\theta)}{|\cos\theta|}+|\cos\theta|\,\hat K_2(\sin\theta)\right)
    \int_{S^1}d\varphi\;|\eta(\theta-\varphi)|^2\\
  &\quad+\frac{1}{2\pi}\,e^{-im\pi}\int_{-\pi}^{\pi}d\theta\;e^{i(n+m)\theta}
    \left(\frac{\hat K_1(\sin\theta)}{|\cos\theta|}-|\cos\theta|\,\hat K_2(\sin\theta)\right)
    \int_{S^1}d\varphi\;\conj{\eta(\theta-\varphi)}\,\eta(\pi-\theta-\varphi)
\end{align*}
where it has bee used that $\cos(\pi-\theta)=-\cos\theta$, \ $\sin(\pi-\theta)=\sin\theta$.

Since
\begin{align*}
  \hat K_1(x)=\sqrt{1-x^2}\quad&\Rightarrow\quad
  \hat K_1(\sin\theta)=\sqrt{1-\sin^2\theta}=|\cos\theta|\\
  \hat K_2(x)=\frac{1}{\sqrt{1-x^2}}\quad&\Rightarrow\quad
  \hat K_2(\sin\theta)=\frac{1}{\sqrt{1-\sin^2\theta}}=\frac{1}{|\cos\theta|}
\end{align*}
\[
  \Rightarrow\quad (\tilde{A}_{\eta})_{m,n}=\frac{1}{2\pi}\int_{-\pi}^{\pi}d\theta\;e^{i(n-m)\theta}
  \int_{S^1}d\varphi\;|\eta(\theta-\varphi)|^2
\]

Expanding: $\eta(\theta)=\frac{1}{\sqrt{2\pi}}\sum_{n\in\Z}a_n\,e^{in\theta}$
\[
  \Rightarrow\quad\int_{S^1}d\varphi\;|\eta(\theta-\varphi)|^2
  =\frac{1}{2\pi}\sum_{n,m\in\Z}\conj{a_n}\,a_m\,e^{i(m-n)\theta}\int_{S^1}d\varphi\;
  e^{-i(m-n)\varphi}
  =\sum_{n\in\Z}|a_n|^2=\,\|\eta\|^2
\]
\[
  \Rightarrow\quad (\tilde{A}_{\eta})_{m,n}=\frac{1}{2\pi}\int_{-\pi}^{\pi}d\theta\;
  e^{i(n-m)\theta}\,\,\|\eta\|^2=\,\|\eta\|^2\,\delta_{m,n}
\]
\[
  \Rightarrow\quad A=\,\|\eta\|^2\,I
\]

\bibliographystyle{ieeetr}

\bibliography{/home/guerrero/MEGA/Universidad/TRABAJO/Bibliography/BibliographyCS.bib}

\end{document}